\documentclass{article}
\usepackage{amsmath}
\usepackage{amssymb}
\usepackage{amscd}

\def\R {\mathbb R}
\def\Z {\mathbb Z}
\def\qed{\hfill$\square$}

\renewcommand{\div}{\mathop{\hbox{\rm div}}}
\newcommand{\grad}{\mathop{\hbox{\rm grad}}}
\newcommand{\curl}{\mathop{\hbox{\rm curl}}}

\newtheorem{theo}{Theorem}
\newtheorem{prop}{Proposition}

\date{}

\begin{document}

\title{Simple finite element schemes\\for the solution of the curl--div system}

\author{Ana Alonso Rodr\'{\i}guez\thanks{Dipartimento di Matematica, Universit\`a di Trento, Trento (Italy).} \and Enrico Bertolazzi\thanks{Dipartimento di Ingegneria Industriale, Universit\`a di Trento, Trento (Italy).}  \and Alberto Valli{$^\ast$}}

\maketitle

\begin{abstract} 
New variational formulations are devised for the curl--div system, and the corresponding finite element approximations are shown to converge. Curl--free and divergence--free finite elements are employed for discretizing the problem.
\end{abstract} 

\section{Introduction, notation and preliminary results}\label{Ch:int}
The curl--div system often appears in electromagnetism (electrostatics, magnetostatics) and in fluid dynamics (rotational incompressible flows, velocity--vorticity formulations). It reads
\begin{equation}\label{eq:cd1}
\left\{\begin{array}{ll}
\curl {\bf u} = {\bf J} &\hbox{\rm in} \ \Omega \\
\div {\bf u}  = g &\hbox{\rm in} \ \Omega\\
{\bf u} \times {\bf n} = {\bf a}  \ \ (\hbox{\rm or} \ \ {\bf u} \cdot {\bf n} = {b}) &\hbox{\rm on} \ \partial \Omega \, ,
\end{array}
\right.
\end{equation}
with in addition some topological conditions assuring uniqueness.

Aim of this paper is to devise some new variational formulations of this problem, which lead to simple and efficient finite element schemes for its numerical approximation.

The main novelty resides in the functional framework we adopt: we look for the solution in the spaces of curl-free or divergence-free vector fields. For the sake of implementation, we also describe in detail how to construct a simple finite element basis for these vector spaces; convergence of the finite element approximations is then easily shown.

\medskip

Let us start with some notation. Let $\Omega$ be a bounded polyhedral domain of $\R^3$ with Lipschitz boundary and let $(\partial \Omega)_0,\dots, (\partial \Omega)_p$ be the connected components of $\partial \Omega$, $(\partial \Omega)_0$ being the external one. We consider a tetrahedral triangulation $\mathcal T_h=(V,E,F,T)$ of $\overline \Omega$, denoting by $V$ the set of vertices, $E$ the set of edges, $F$ the set of faces and $T$ the set of tetrahedra of $\mathcal T_h$.

We will use these spaces of finite elements (see Monk~\cite{Monk} for a complete presentation): the space $L_h$ of continuous piecewise-linear elements, with dimension $n_v$, the number of vertices in $\mathcal T_h$; the space $N_h$ of N\'ed\'elec edge elements of degree 1, with dimension $n_e$, the number of edges in $\mathcal T_h$; the space $RT_h$ of Raviart-Thomas elements of degree 1, with dimension $n_f$, the number of faces in $\mathcal T_h$; the space $PC_h$ of piecewise-constant elements, with dimension $n_t$, the number of tetrahedra in $\mathcal T_h$.

The following inclusions are well-known: 
$$
\begin{array}{ll}
L_h \subset H^1(\Omega) \ \ , \ \ N_h \subset H(\curl; \Omega) \ \ , \ \ RT_h \subset H(\div; \Omega) \ \, \ \ PC_h \subset L^2(\Omega) \, ,
\end{array}
$$
where
$$
\begin{array}{l}
H^1(\Omega)=\{ \phi \in L^2(\Omega) \, | \, \grad \phi \in (L^2(\Omega))^3 \} \, , \\[.2cm]
H(\curl; \Omega)=\{ {\bf v} \in (L^2(\Omega))^3 \, | \, \curl {\bf v} \in (L^2(\Omega))^3 \} \, ,\\[.2cm]
H(\div; \Omega)=\{ {\bf v} \in (L^2(\Omega))^3 \, | \, \div {\bf v} \in L^2(\Omega) \} \, .
\end{array}
$$
Moreover $\grad L_h \subset N_h$, $\curl N_h \subset RT_h$, and $\div RT_h \subset PC_h$.
The basis of $L_h$ is denoted by $\{ \psi_{h,1}, \dots, \psi_{h,n_v} \}$, with 
${\psi}_{h,i} (v_j) = \delta_{i,j}$
for $1 \le i,j \le n_v$; the basis of $N_h$ is denoted by $\{ {\bf w}_{h,1}, \dots, {\bf w}_{h,n_e} \}$, with 
$\int_{e_j} {\bf w}_{h,i} \cdot \boldsymbol \tau = \delta_{i,j}$
for $1 \le i,j \le n_e$; the basis of  $RT_h$ is denoted by $\{ {\bf r}_{h,1}, \dots, {\bf r}_{h,n_f} \}$, with 
$\int_{f_m} {\bf r}_{h,l} \cdot \boldsymbol \nu = \delta_{l,m}$
for $1 \le l,m \le n_f$.

Fixing a total ordering $v_1,\dots,v_{n_v}$ of the elements of $V$, an orientation on the elements of $E$ and $F$ is induced: if the end points of $e_j$ are $v_a$ and $v_b$ for some $a, b \in \{ 1,\dots, n_v \}$ with $a < b$, then the oriented edge $e_j$ will be denoted by $[v_a,v_b]$, with unit  tangent vector $\boldsymbol \tau= \frac{v_b-v_a}{|v_b-v_a|}$; if the face $f_m$ has vertices $v_a$, $v_b$ and $v_c$ with $a < b <c$, the oriented face $f_m$ will be denoted by $[v_a,v_b,v_c]$ and its unit normal vector $\boldsymbol \nu = \frac{(v_b - v_a) \times (v_c-v_a)}{|(v_b - v_a) \times (v_c-v_a)|}$ is obtained by the right hand rule. 

We finally introduce a set of 1-cycles in ${\mathcal T}_h$, denoted by $\{ \sigma_n \}_{n=1}^{g}$, that are representatives of a basis of the first homology group ${\mathcal H}_1(\overline\Omega,\Z)$ (whose rank is therefore equal to $g$): in other words, this set is a maximal set of non-bounding 1-cycles in ${\mathcal T}_h$. Let us recall that an explicit and efficient construction of the 1-cycles $\{ \sigma_n \}_{n=1}^{g}$ is given by Hiptmair and Ostrowski~\cite{HO02}. For a more detailed presentation of the homological concepts that are useful in the numerical approximation of PDEs, see, e.g., Bossavit~\cite{Boss98}, Hiptmair~\cite{Hip02}, Gross and Kotiuga~\cite{GK04}; see also Benedetti, Frigerio and Ghiloni~\cite{BFG10}, Alonso Rodr{\'\i}guez et al.~\cite{ABGV13}.

\section{The curl--div system with assigned tangential component on the boundary}

Let $\boldsymbol \eta$ be a symmetric matrix, uniformly positive definite in $\Omega$, with entries belonging to $L^\infty(\Omega)$. Given ${\bf J} \in (L^2(\Omega))^3$, $g \in L^2(\Omega)$, ${\bf a} \in H^{-1/2}(\div_\tau;\partial \Omega)$, $\boldsymbol\alpha \in \R^p$, we look for ${\bf u} \in (L^2(\Omega))^3$ such that  
\begin{equation}\label{eq:cd_tan}
\left\{\begin{array}{ll}
\curl (\boldsymbol \eta {\bf u}) = {\bf J}  &\hbox{\rm in} \ \Omega\\
\div {\bf u}  = g &\hbox{\rm in} \ \Omega\\
(\boldsymbol \eta {\bf u}) \times {\bf n} = {\bf a} &\hbox{\rm on} \ \partial\Omega \\
\int_{(\partial \Omega)_r} {\bf u} \cdot {\bf n} =\alpha_r &\hbox{\rm for each} \ r = 1,\dots,p \, .
\end{array}
\right.
\end{equation}
The data satisfy the necessary conditions  $\div {\bf J} = 0$ in $\Omega$, $\int_{\Omega} {\bf J} \cdot \boldsymbol {\rho} + \int_{\partial \Omega} {\bf a} \cdot \boldsymbol\rho = 0$ for each $\boldsymbol\rho \in {\mathcal H}(m)$,  and ${\bf J} \cdot {\bf n}= \div_\tau {\bf a}$ on $\partial \Omega$, where $ {\mathcal H}(m)$ is the space of Neumann harmonic fields, namely,
$$
 {\mathcal H}(m) = \{\boldsymbol\rho \in (L^2(\Omega))^3 \, | \, \curl \boldsymbol\rho = {\bf 0} \ \hbox{\rm in} \ \Omega, \div \boldsymbol\rho = 0 \ \hbox{\rm in} \ \Omega,  \boldsymbol\rho \cdot  {\bf n}= 0 \ \hbox{\rm on} \ \partial \Omega\} \, .
$$

By means of a variational approach, Saranen~\cite{Sa82}, \cite{Sa82b} has shown that this problem has a unique solution (see also the results proved in Alonso Rodr{\'\i}guez and Valli~\cite[Sect.\ A3]{AVbook}). The method is based on a splitting of the solution in three terms, writing 
$$
\boldsymbol\eta {\bf u} = \boldsymbol\eta \curl {\bf q} +  \grad \chi + {\bf h} \, .
$$
Here ${\bf h}$ is a generalized harmonic field satisfying $\curl {\bf h} = {\bf 0}$ in $\Omega$, $\div(\boldsymbol\eta^{-1} {\bf h})=0$ in $\Omega$ and ${\bf h} \times {\bf n} = {\bf 0}$ on $\partial \Omega$, namely, it is an element of a finite dimensional vector space of dimension $p$; ${\bf q}$ is a solution to $\curl(\boldsymbol\eta \curl {\bf q})= {\bf J}$ in $\Omega$ and $(\boldsymbol\eta \curl {\bf q}) \times {\bf n} = {\bf a}$ on $\partial \Omega$; $\chi$ is the solution to $\div(\boldsymbol\eta^{-1} \grad \chi) = g$ in $\Omega$ and $\chi=0$ on $\partial \Omega$.

Since a solution ${\bf q}$ to $\curl(\boldsymbol\eta \curl {\bf q})= {\bf J}$ in $\Omega$ and $(\boldsymbol\eta \curl {\bf q}) \times {\bf n} = {\bf a}$ on $\partial \Omega$ is not unique (${\bf q} + \grad \phi$ is still a solution), other equations have to be added. Typically, one imposes the gauging conditions $\div {\bf q}=0$ in $\Omega$, ${\bf q}\cdot {\bf n}=0$ on $\partial \Omega$ and ${\bf q} \bot {\mathcal H}(m)$. 

This leads to two variational problems: a standard Dirichlet boundary value problem for $\chi$, and a constrained problem for ${\bf q}$  (the determination of the harmonic field ${\bf h}$ also needs some additional work, but it is an easy finite dimensional problem).
Before proceeding, let us look at some possible numerical approaches for approximating these two problems.

The first one is a standard elliptic problem. Numerical approximation can be performed by scalar nodal elements in $H^1(\Omega)$, looking for the unknown $\chi$ and then computing its gradient, or by means of a (more expensive) mixed method in $H(\div;\Omega) \times L^2(\Omega)$, in which $\grad \chi \in H(\div;\Omega)$ is directly computed as an auxiliary unknown.

Concerning the problem related to the vector field ${\bf q}$, a first choice is to work in $H(\curl;\Omega) \cap H(\div;\Omega)$, hence with globally-continuous nodal finite elements for each component of ${\bf q}$; the drawback is that, in the presence of re-entrant corners, the solution is singular (it does not belong to $(H^1(\Omega))^3$) and $(H^1(\Omega))^3$ is a closed subspace of $H(\curl;\Omega) \cap H(\div;\Omega)$, hence in that case a finite element scheme cannot be convergent (see, e.g. Costabel, Dauge and Nicaise~\cite{CDN04}).

An alternative method is to formulate the problem as a saddle-point problem in $H(\curl;\Omega)$, in which the divergence constraint is imposed in a week sense, introducing a scalar Lagrange multiplier; in this way the number of degrees of freedom is rather high, as, besides an edge approximation of the vector field ${\bf q}$, one has also to consider a nodal approximation of the scalar Lagrange multiplier. Moreover, the resulting algebraic problem is associated to a non-definite matrix.

A way for avoiding the introduction of a Lagrange multiplier is to solve the equation $\curl(\boldsymbol\eta \curl {\bf q})= {\bf J}$ in $\Omega$ by using edge elements without any gauge; though the conjugate gradient method has shown to be a viable tool for solving the associated algebraic problem (see, e.g., B{\'\i}r\'o~\cite{BIR99}), the price to pay is a worst convergence rate, since the matrix to deal with is singular.

\subsection{A new variational formulation}\label{varfor}

The previous  discussion should have explained why our aim here is to find a different variational formulation for problem (\ref{eq:cd_tan}), a formulation that will be more suitable for numerical approximation. Our method is somehow related to the so-called tree--cotree gauge used for the numerical approximation of eddy current  and magnetostatic problems (see, e.g, Ren and Razek~\cite{RR90}, Albanese and Rubinacci~\cite{AR98}).

The first step of the procedure is to find a vector field ${\bf u}^\star \in (L^2(\Omega))^3$ satisfying
\begin{equation}\label{eq:u*}
\left\{\begin{array}{ll}
\div {\bf u}^\star  = g &\hbox{\rm in} \ \Omega\\
\int_{(\partial \Omega)_r} {\bf u}^\star \cdot {\bf n} =\alpha_r &\hbox{\rm for each} \ r = 1,\dots,p \, .
\end{array}
\right.
\end{equation}
Such a vector field does exist: for instance, one can think to take ${\bf J} = {\bf 0}$ and ${\bf a} = {\bf 0}$ in (\ref{eq:cd_tan}), or any choice of ${\bf J}$ and ${\bf a}$ satisfying the compatibility conditions (indeed, we will not assume in the sequel that $\curl (\boldsymbol \eta {\bf u}^\star) = {\bf 0}$ or $(\boldsymbol \eta {\bf u}^\star) \times {\bf n} = {\bf 0}$).

The vector field ${\bf W}= {\bf u} - {\bf u}^\star$ satisfies
\begin{equation}\label{eq:cd_tan2}
\left\{\begin{array}{ll}
\curl (\boldsymbol \eta {\bf W}) = {\bf J} - \curl (\boldsymbol \eta {\bf u}^\star)  &\hbox{\rm in} \ \Omega\\
\div {\bf W}  = 0 &\hbox{\rm in} \ \Omega\\
(\boldsymbol \eta {\bf W}) \times  {\bf n} = {\bf a} - (\boldsymbol \eta {\bf u}^\star) \times  {\bf n} &\hbox{\rm on} \ \partial\Omega \\
\int_{(\partial \Omega)_r} {\bf W} \cdot {\bf n} =0 &\hbox{\rm for each} \ r = 1,\dots,p \, ,
\end{array}
\right.
\end{equation}
and the second step of the procedure is finding a simple variational formulation of this problem.

Multiplying the first equation by a test function ${\bf v} \in H(\curl;\Omega)$, integrating in $\Omega$ and integrating by parts we find: 
$$
\begin{array}{ll}
\int_\Omega {\bf J} \cdot {\bf v} \!\!\!&= \int_\Omega \curl [\boldsymbol \eta ({\bf W} + {\bf u}^\star)] \cdot {\bf v} \\[1ex]
&= \int_\Omega\boldsymbol \eta ({\bf W} + {\bf u}^\star) \cdot \curl {\bf v} - \int_{\partial \Omega} 
 [\boldsymbol \eta ({\bf W} + {\bf u}^\star) \times {\bf n}] \cdot {\bf v}
\\[1ex]
&= \int_\Omega\boldsymbol \eta {\bf W}  \cdot \curl {\bf v} + \int_\Omega\boldsymbol \eta {\bf u}^\star \cdot \curl {\bf v} - \int_{\partial \Omega} 
 {\bf a} \cdot {\bf v} \, .
\end{array}
$$
Let us introduce the space
\begin{equation}\label{eq:space}
\begin{array}{ll}
{\mathcal W}_0 = \{{\bf v} \in H(\div;\Omega) \, | \, \div{\bf v} = 0 \ \hbox{\rm in} \ \Omega, \\
\hspace{3.5cm} \int_{(\partial \Omega)_r} {\bf v} \cdot {\bf n} =0  \ \hbox{\rm for each} \ r = 1,\dots,p\} \, .
\end{array}
\end{equation}
Note that this space can be written as ${\mathcal W}_0 = \curl [H(\curl;\Omega)]$: in fact, the inclusion $\curl [H(\curl;\Omega)] \subset {\mathcal W}_0$ is obvious, while the inclusion ${\mathcal W}_0 \subset \curl [H(\curl;\Omega)]$ is the classical result concerning vector potentials (see, e.g., Cantarella, DeTurck and Gluck~\cite{CDG02}, Alonso Rodr{\'\i}guez and Valli~\cite{AV14}).
The vector field ${\bf W}$ is thus a solution to 
\begin{equation}\label{eq:varpb}
\begin{array}{ll}
{\bf W} \in {\mathcal W}_0 \ : \ \int_\Omega\boldsymbol \eta {\bf W}  \cdot \curl {\bf v} = \int_\Omega {\bf J} \cdot {\bf v}- \int_\Omega\boldsymbol \eta {\bf u}^\star \cdot \curl {\bf v}\\[1ex]
\hspace{6cm} + \int_{\partial \Omega} {\bf a} \cdot {\bf v}  \ \ \ \ \forall \ {\bf v} \in H(\curl;\Omega) \, .
\end{array}
\end{equation}
More precisely, ${\bf W}$ is the {\sl unique} solution of that problem: in fact, assuming ${\bf J}= {\bf u}^\star= {\bf a}={\bf 0}$, and taking ${\bf v}$ such that $\curl {\bf v} = {\bf W}$, it follows at once $\int_\Omega \boldsymbol \eta {\bf W} \cdot {\bf W} = 0$, hence ${\bf W}= {\bf 0}$.

\subsection{Finite element approximation}\label{FEMapprox}

The finite element approximation of problem (\ref{eq:cd_tan}) can be performed in two steps, too. The first one is finding a finite element potential ${\bf u}^\star_h \in RT_h$ such that
\begin{equation}\label{eq:u*h}
\left\{\begin{array}{ll}
\div {\bf u}^\star_h  = g_h &\hbox{\rm in} \ \Omega\\
\int_{(\partial \Omega)_r} {\bf u}^\star_h \cdot {\bf n} =\alpha_r &\hbox{\rm for each} \ r = 1,\dots,p \, ,
\end{array}
\right.
\end{equation}
where $g_h \in PC_h$ is the piecewise-constant interpolant $I^{PC}_h g$ of $g$. This can be done by means of a simple and efficient algorithm as shown in Alonso Rodr{\'\i}guez and Valli~\cite{AV14}.

The second step concerns the numerical approximation of problem (\ref{eq:varpb}). Here the main issue is to determine a finite element subspace of  ${\mathcal W}_{0}$, and a suitable finite element basis. The natural choice is clearly
\begin{equation}\label{eq:space_h}
\begin{array}{ll}
{\mathcal W}_{0,h} = \{{\bf v}_h \in RT_h \, | \, \div{\bf v}_h = 0 \ \hbox{\rm in} \ \Omega, \\
\hspace{3.5cm} \int_{(\partial \Omega)_r} {\bf v}_h \cdot {\bf n} =0  \ \hbox{\rm for each} \ r = 1,\dots,p\} \, .
\end{array}
\end{equation}
For the ease of notation let us set $n_Q=n_e - (n_v - 1)$. As proved in Alonso Rodr{\'\i}guez et al.~\cite{ACGV15}, the dimension of ${\mathcal W}_{0,h}$ is equal to $n_Q - g$, and a basis is given by the curls of suitable N\'ed\'elec elements belonging to $N_h$. 

To make clear this point, following Alonso Rodr{\'\i}guez et al.~\cite{ACGV15}, some notations are necessary.
As shown in Hiptmair and Ostrowski~\cite{HO02} (see also Alonso Rodr{\'\i}guez et al.~\cite{ABGV13}), it is possible to construct a set of 1-cycles $\{\sigma_n\}_{n=1}^g$, representing a basis of the first homology group ${\mathcal H}_1(\overline\Omega,\Z)$, as
a formal sum of edges in $\mathcal T_h$ with integer coefficients. More precisely, let us consider the graph given by the vertices and the edges of $\mathcal T_h$ on $\partial \Omega$. The number of connected components of this graph coincides with the number of connected components of $\partial \Omega$. For each $r=0,1,\dots,p$ let $S_{\partial \Omega}^r=(V_{\partial \Omega}^r,M_{\partial \Omega}^r)$ be a spanning tree of the corresponding connected component of the graph. Then consider the graph $(V,E)$, given by all the vertices and edges of $\mathcal T_h$, and a spanning tree $S=(V,M)$ of this graph such that $M_{\partial \Omega}^r \subset M$ for each $r=0,1,\dots,p$. Let us order the edges in such a way that the edge $e_l$ belongs to the cotree of $S$ for $l=1,\dots,n_Q$ and the edge $e_{n_Q+i}$ belongs to the tree $S$ for $i=1,\dots,n_v - 1$. In particular, denote by $e_q$, $q=1,\ldots,2g$, the set of edges of $\partial \Omega$, constructed by Hiptmair and Ostrowski~\cite{HO02}, that have the following properties: they all belong to the cotree, and each one of them ``closes" a 1-cycle $\gamma_q$ that is a representative of a basis of the first homology group ${\mathcal H}_1(\partial\Omega,\Z)$ (whose rank is indeed equal to $2g$).  
With this notation, we recall that the 1-cycles $\sigma_n$ can be expressed as the formal sum
\begin{equation}\label{A}
\sigma_n = \sum_{q =1}^{2g} A_{n,q} \gamma_q= \sum_{q =1}^{2g} A_{n,q} e_q + \sum_{i=n_Q+1}^{n_e} a_{n,i} e_i \, ,
\end{equation}
for suitable and explicitly computable integers $A_{n,q}$. 

The idea that leads to the construction of the basis of ${\mathcal W}_{0,h}$ is now the following: first, consider the set 
$$
\{\curl {\bf w}_{h,l}\}_{l=2g+1}^{n_Q} \, ,
$$
Then look for $g$ functions ${\bf z}_{h,\lambda} \in RT_h$, $\lambda=1,\dots,g$, of the form
$$
{\bf z}_{h,\lambda}=\sum_{\upsilon=1}^{2g} c^{(\lambda)}_{\upsilon} \curl {\bf w}_{h,\upsilon} \, ,
$$
where the linearly independent vectors ${\bf c}^{(\lambda)} \in \R^{2g}$ are chosen in such a  way that
$$
\oint_{\sigma_n}\left(\sum_{\upsilon=1}^{2g} c^{(\lambda)}_{\upsilon} {\bf w}_{h,\upsilon} \right) \cdot d{\bf s} =0
$$
for $n=1,\dots,g$. This can be done since $\sigma_n$ is formed by the ``closing" edges $e_q$, $q=1,\ldots,2g$, and by edges belonging to the spanning tree, so that
$$
\begin{array}{ll}
\displaystyle \oint_{\sigma_n}\left( \sum_{\upsilon=1}^{2g}c^{(\lambda)}_{\upsilon} {\bf w}_{h,\upsilon}  \right) \cdot d{\bf s} \!\!\!&=
\displaystyle \sum_{q =1}^{2g} A_{n,q} \int_{e_q}\left( \sum_{\upsilon=1}^{2g}c^{(\lambda)}_{\upsilon} {\bf w}_{h,\upsilon} \right) \cdot \boldsymbol \tau=\displaystyle\sum_{q =1}^{2g} A_{n,q}  c^{(\lambda)}_{q} \, ,
\end{array}
$$
and the matrix $A \in \Z^{g \times 2g}$ with entries $A_{n,q}$ has rank $g$ (see Hiptmair and Ostrowski~\cite{HO02}, Alonso Rodr{\'\i}guez et al.~\cite[Sect.\ 6]{ABGV13}). Thus we only have to determine a basis ${\bf c}^{(\lambda)} \in \R^{2g}$ of the kernel of $A$, $\lambda=1,\ldots,g$. An easy way for determining these vectors ${\bf c}^{(\lambda)}$ is presented in Alonso Rodr{\'\i}guez et al.~\cite{ACGV15}. 

\begin{prop}\label{theo:basis}
The vector fields
$$
\{\curl {\bf w}_{h,l}\}_{l=2g+1}^{n_Q}  \cup \Big\{\curl \Big(\sum_{\upsilon=1}^{2g} c^{(\lambda)}_{\upsilon} {\bf w}_{h,\upsilon}\Big)\Big\}_{\lambda=1}^g \subset {\mathcal W}_{0,h}
$$
are linearly independent. 
\end{prop}
{\bf Proof.}
The proof is in Alonso Rodr{\'\i}guez et al.~\cite[Prop.\ 2]{ACGV15}.  \qed  

\vspace{.12cm}

The straightforward consequence of this result is that, since the number of the vector fields in Proposition~\ref{theo:basis} is $n_Q - g$, the dimension of ${\mathcal W}_{0,h}$, they are a basis of ${\mathcal W}_{0,h}$. Let us denote this basis by $\{\curl \boldsymbol\omega_{h,l}\}_{l=g+1}^{n_Q}$, with 
\begin{equation}\label{basis}
\boldsymbol\omega_{h,l}= \left\{\begin{array}{ll}{\bf w}_{h,l} &\hbox{\rm for} \ l=2g+1,\ldots,n_Q \\ [1ex]
\sum_{\upsilon=1}^{2g} c^{(l-g)}_{\upsilon} {\bf w}_{h,\upsilon} &\hbox{\rm for} \ l=g+1,\ldots,2g \, .\end{array}
 \right.
 \end{equation}

\begin{prop}
The vector fields $\{\boldsymbol\omega_{h,l}\}_{l=g+1}^{n_Q}$ are linearly independent.
\end{prop}
{\bf Proof.} Suppose we have $\sum_{l=g+1}^{n_Q} \kappa_l \boldsymbol\omega_{h,l} = {\bf 0}$ for some $\kappa_l$. This can be rewritten as
$$
\begin{array}{ll}
{\bf 0}\!\!\!&= \sum_{l=2g+1}^{n_Q} \kappa_l {\bf w}_{h,l} + \sum_{l=g+1}^{2g} \kappa_l\Big(\sum_{\upsilon=1}^{2g} c^{(l-g)}_{\upsilon} {\bf w}_{h,\upsilon}\Big) \\[1ex]
&= \sum_{l=2g+1}^{n_Q} \kappa_l {\bf w}_{h,l} + \sum_{\upsilon=1}^{2g}\Big(\sum_{l=g+1}^{2g} \kappa_l c^{(l-g)}_{\upsilon} \Big){\bf w}_{h,\upsilon} \, ,
\end{array}
$$
thus $\kappa_l =0$ for $l=2g+1,\ldots,n_Q$ and $\sum_{l=g+1}^{2g} \kappa_l c^{(l-g)}_{\upsilon} =0$ for $\upsilon=1,\ldots2g$, as $\{{\bf w}_{h,l}\}_{l=1}^{n_Q}$ are linearly independent. Since the vectors ${\bf c}^{(l-g)} \in \R^{2g}$, $l=g+1,\ldots,2g$, are linearly independent, we also obtain $\kappa_l=0$ for $l=g+1,\ldots,2g$, and the result follows. \qed

\vspace{0.12cm}

We are now in a position to formulate the finite element approximation of (\ref{eq:varpb}), that reads as follows:
\begin{equation}\label{eq:varpb_fin}
\begin{array}{ll}
{\bf W}_h \in {\mathcal W}_{0,h} \ : \ \int_\Omega\boldsymbol \eta {\bf W}_h \cdot \curl {\bf v}_h = \int_\Omega {\bf J} \cdot {\bf v}_h- \int_\Omega\boldsymbol \eta {\bf u}^\star_h \cdot \curl {\bf v}_h\\[1ex]
\hspace{6cm} + \int_{\partial \Omega} {\bf a} \cdot {\bf v}_h  \ \ \ \ \forall \ {\bf v}_h \in N_h^\star\, ,
\end{array}
\end{equation}
where 
\begin{equation}\label{FEspace_red}
N^\star_h = \hbox{\rm span} \{\boldsymbol\omega_{h,l}\}_{l=g+1}^{n_Q} \, .
\end{equation}
The corresponding algebraic problem is a square linear system of dimension $n_Q - g$, and it is uniquely solvable. In fact,  we note that ${\mathcal W}_{0,h} = \curl  N_h^\star$, hence we can choose ${\bf v}_h^\star \in N_h^\star$ such that $\curl {\bf v}_h^\star = {\bf W}_h$; from $(\ref{eq:varpb_fin})$ we find at once ${\bf W}_h = {\bf 0}$, provided that ${\bf J}= {\bf u}^\star_h= {\bf a}={\bf 0}$.

The convergence of this finite element scheme is easily shown by standard arguments. For the ease of reading, let us present the proof.

\begin{theo} Let ${\bf W} \in {\mathcal W}_0$ and ${\bf W}_h\in {\mathcal W}_{0,h}$ be the solutions of problem (\ref{eq:varpb}) and (\ref{eq:varpb_fin}), respectively. Set ${\bf u} = {\bf W} + {\bf u}^\star$ and ${\bf u}_h= {\bf W}_h + {\bf u}_h^\star$, where ${\bf u}^\star \in H(\div;\Omega)$ and ${\bf u}_h^\star \in RT_h$ are solutions to problem (\ref{eq:u*}) and (\ref{eq:u*h}), respectively. Assume that ${\bf u}$ is regular enough, so that the interpolant $I^{RT}_h {\bf u}$ is defined. Then the following error estimate holds
\begin{equation}\label{eq:error}
\|{\bf u}-{\bf u}_h\|_{H(\scriptsize \div;\Omega)} \le c_0 (\|{\bf u} - I^{RT}_h {\bf u}\|_{L^2(\Omega)} + \|g - I^{PC}_h g\|_{L^2(\Omega)}) \, .
\end{equation}
\end{theo}
{\bf Proof.} Since $N_h^\star \subset H(\curl;\Omega)$, we can choose ${\bf v}= {\bf v}_h \in N_h^\star$ in (\ref{eq:varpb}). By subtracting (\ref{eq:varpb_fin}) from (\ref{eq:varpb}) we end up with
$$
\int_\Omega \boldsymbol\eta [({\bf W} + {\bf u}^\star) - ({\bf W}_h + {\bf u}_h^\star)] \cdot \curl {\bf v}_h = 0 \qquad \forall \ {\bf v}_h \in N_h^\star \, ,
$$
namely, 
\begin{equation}\label{eq:consist}
\int_\Omega \boldsymbol\eta ({\bf u} -  {\bf u}_h)\cdot \curl {\bf v}_h = 0 \qquad \forall \ {\bf v}_h \in N_h^\star \, .
\end{equation}
Then, recalling that ${\mathcal W}_{0,h} = \curl N_h^\star$, therefore ${\bf W}_h = \curl {\bf v}_h^\star$ for a suitable ${\bf v}^\star_h \in N^\star_h$, and using (\ref{eq:consist}) we find
$$
\begin{array}{ll}
c_1 \|{\bf u}-{\bf u}_h\|_{L^2(\Omega)}^2 \!\!\!&\le \int_\Omega \boldsymbol \eta ({\bf u}-{\bf u}_h) \cdot ({\bf u}-{\bf u}_h)\\ [1ex]
&=  \int_\Omega \boldsymbol \eta ({\bf u}-{\bf u}_h) \cdot ({\bf u} - {\bf W}_h - {\bf u}_h^\star) \\[1ex]
&=  \int_\Omega \boldsymbol \eta ({\bf u}-{\bf u}_h) \cdot ({\bf u} - \curl{\bf v}_h^\star - {\bf u}_h^\star) \\ [1ex]
&=  \int_\Omega \boldsymbol \eta ({\bf u}-{\bf u}_h) \cdot ({\bf u} - \curl{\bf v}_h - {\bf u}_h^\star) \\[1ex]
&\le c_2 \|{\bf u}-{\bf u}_h\|_{L^2(\Omega)} \|{\bf u} - \boldsymbol\Phi_h  - {\bf u}_h^\star\|_{L^2(\Omega)} \quad \forall \ {\boldsymbol\Phi}_h \in {\mathcal W}_{0,h} \, .
\end{array}
$$
We can choose ${\boldsymbol\Phi}_h = (I^{RT}_h {\bf u}- {\bf u}_h^\star) \in {\mathcal W}_{0,h}$; in fact, $\div (I^{RT}_h {\bf u}) = I^{PC}_h(\div {\bf u}) = I^{PC}_h g=g_h$ and $\int_{(\partial \Omega)_r} I^{RT}_h {\bf u} \cdot {\bf n} = \int_{(\partial \Omega)_r}  {\bf u} \cdot {\bf n} = \alpha_r$ for each $r=1,\ldots,p$.
Then  it follows at once $\|{\bf u}-{\bf u}_h\|_{L^2(\Omega)} \le \frac{c_2}{c_1} \|{\bf u} - I^{RT}_h {\bf u}\|_{L^2(\Omega)}$.
Moreover, $\div({\bf u}-{\bf u}_h) = g - g_h = g-I^{PC}_h g$, and the thesis is proved. \qed

\vspace{0.12cm}

Note that a sufficient condition for defining the interpolant of ${\bf u}$ is that ${\bf u} \in (H^{\frac{1}{2} + \delta}(\Omega))^3$, $\delta >0$ (see Monk~\cite[Lemma 5.15]{Monk}). This is satisfied if, e.g., $\boldsymbol \eta$ is a scalar Lipschitz function in $\overline \Omega$ and ${\bf a} \in (H^\gamma(\partial \Omega))^3$, $\gamma > 0$ (see Alonso and Valli~\cite{AV99}).

\subsection{The algebraic problem}

The solution ${\bf W}_h \in {\mathcal W}_{0,h}$ can be written in terms of the basis as ${\bf W}_h = \sum_{l=g+1}^{n_Q} W_l \curl \boldsymbol\omega_{h,l}$. Hence the finite dimensional problem (\ref{eq:varpb_fin}) can be rewritten as
\begin{equation}\label{eq:matrix}
\begin{array}{ll}
\displaystyle\sum_{l=g+1}^{n_Q} W_l \int_\Omega \boldsymbol{\eta} \curl \boldsymbol\omega_{h,l} \cdot \curl \boldsymbol\omega_{h,m} = \int_\Omega {\bf J} \cdot \boldsymbol\omega_{h,m}- \int_\Omega\boldsymbol \eta {\bf u}^\star_h \cdot \curl \boldsymbol\omega_{h,m}\\
\hspace{6cm} + \displaystyle\int_{\partial \Omega} {\bf a} \cdot \boldsymbol\omega_{h,m}   \, ,
\end{array}
\end{equation}
for each  $m=g+1,\ldots,n_Q$. 

\begin{theo}
The matrix ${\bf K}^\star$ with entries 
$$
K_{ml}^\star = \int_\Omega \boldsymbol{\eta} \curl \boldsymbol\omega_{h,l} \cdot \curl \boldsymbol\omega_{h,m}
$$
is symmetric and positive definite.
\end{theo}
{\bf Proof.} It is enough to recall that the vector fields $\{ \curl \boldsymbol\omega_{h,l}\}_{l=g+1}^{n_Q}$ are linearly independent (see Proposition~\ref{theo:basis}). More precisely, they are a basis of ${\mathcal W}_{0,h}$, hence ${\bf Q}^\star$ is the mass matrix in ${\mathcal W}_{0,h}$ with weight $\boldsymbol{\eta}$. \qed

\section{The curl--div system with assigned normal component on the boundary}

Let $\boldsymbol \mu$ be a symmetric matrix, uniformly positive definite in $\Omega$, with entries belonging to $L^\infty(\Omega)$. Given ${\bf J} \in (L^2(\Omega))^3$, $g \in L^2(\Omega)$, ${b} \in H^{-1/2}(\partial \Omega)$, $\boldsymbol\beta \in \R^g$, we look for ${\bf u} \in (L^2(\Omega))^3$ such that 
\begin{equation}\label{eq:cd_nor}
\left\{\begin{array}{ll}
\curl {\bf u} = {\bf J}  \\
\div (\boldsymbol \mu{\bf u})  = g \\
\boldsymbol \mu{\bf u} \cdot {\bf n} = {b} \\
\oint_{\sigma_n} {\bf u} \cdot d{\bf s} = \beta_n &\hbox{\rm for each} \ n = 1,\dots,g  \, ,
\end{array}
\right.
\end{equation}
where the data satisfy the necessary conditions $\div {\bf J} = 0$ in $\Omega$, $\int_\Omega g= \int_{\partial \Omega} b$; moreover, in order that the line integral of ${\bf u}$ on $\sigma_n$ has a meaning, we also assume that ${\bf J} \cdot {\bf n} = 0$ on $\partial \Omega$ (which is more restrictive than the necessary condition $\int_{(\partial \Omega)_r} {\bf J} \cdot {\bf n} = 0$ for each $r=1,\ldots,p$).

The variational approach proposed by Saranen~\cite{Sa82}, \cite{Sa82b} shows that this problem has a unique solution (see also Alonso Rodr{\'\i}guez and Valli~\cite[Sect.\ A3]{AVbook}). The method is similar to that employed for the problem in which the tangential component of ${\bf u}$ is assigned on the boundary; thus we do not give other details here and refer to the presentation given before.

\subsection{A new variational formulation}\label{varfor2}

The variational formulation of the curl--div system with assigned normal component on the boundary that we present here is similar to the one we have proposed in Alonso Rodr{\'\i}guez et al.~\cite{ABGV13} for the problem of magnetostatics. However, we think it can be interesting for its specific simplicity, as here we will formulate the problem in the space ${\mathcal V}_0 = \grad [H^1(\Omega)]$, while in \cite{ABGV13} it was set in the space $H^0(\curl;\Omega) =\{ {\bf v} \in (L^2(\Omega))^3 \, | \, \curl {\bf v} = {\bf 0} \ \hbox{\rm in} \ \Omega\}$, which in the general topological case is more complicate to approximate.

Also in this case we need a preliminary step: to find a vector field ${\bf u}^{\ast} \in (L^2(\Omega))^3$ satisfying
\begin{equation}\label{eq:u*2}
\left\{\begin{array}{ll}
\curl {\bf u}^\ast  = {\bf J} &\hbox{\rm in} \ \Omega\\
\oint_{\sigma_n} {\bf u}^\ast \cdot d{\bf s} = \beta_n&\hbox{\rm for each} \ n = 1,\dots,g \, .
\end{array}
\right.
\end{equation}
This vector field does exist: for instance, one can choose $g = 0$ and $b= 0$ in (\ref{eq:cd_nor}), or any choice of ${g}$ and ${b}$ satisfying the compatibility condition (indeed, we do not need to assume in the sequel that $\div (\boldsymbol \mu {\bf u}^\ast) = {0}$ or $(\boldsymbol \mu {\bf u}^\ast) \cdot {\bf n} = {0}$).

The vector field ${\bf V}= {\bf u} - {\bf u}^\ast$ satisfies
\begin{equation}\label{eq:cd_nor2}
\left\{\begin{array}{ll}
\curl  {\bf V} = {\bf 0}  &\hbox{\rm in} \ \Omega\\
\div (\boldsymbol \mu{\bf V})  = g - \div (\boldsymbol \mu{\bf u}^\ast) &\hbox{\rm in} \ \Omega\\
(\boldsymbol \mu {\bf V}) \cdot  {\bf n} = {b} - (\boldsymbol \mu {\bf u}^\ast) \cdot  {\bf n} &\hbox{\rm on} \ \partial\Omega \\
\oint_{\sigma_n} {\bf V} \cdot d{\bf s} = 0 &\hbox{\rm for each} \ n = 1,\dots,g \, ,
\end{array}
\right.
\end{equation}
and now we only have to find a variational formulation of this problem.

Multiplying the second equation by a test function ${\varphi} \in H^1(\Omega)$, integrating in $\Omega$ and integrating by parts we find: 
$$
\begin{array}{ll}
\int_\Omega {g} \, {\varphi} \!\!\!&= \int_\Omega \div [\boldsymbol \mu ({\bf V} + {\bf u}^\ast)] \, {\varphi} \\[1ex]
&= - \int_\Omega\boldsymbol \mu ({\bf V} + {\bf u}^\ast) \cdot \grad {\varphi} + \int_{\partial \Omega} 
 [\boldsymbol \mu ({\bf V} + {\bf u}^\ast) \cdot {\bf n}] \, \varphi
\\[1ex]
&= -\int_\Omega\boldsymbol \mu {\bf V}  \cdot \grad {\varphi} - \int_\Omega\boldsymbol \mu {\bf u}^\ast \cdot \grad \varphi + \int_{\partial \Omega} 
 {b} \, {\varphi} \, .
\end{array}
$$
Let us introduce the space
\begin{equation}\label{eq:space2}
\begin{array}{ll}
{\mathcal V}_0 = \{{\bf v} \in H(\curl;\Omega) \, | \, \curl{\bf v} = {\bf 0} \ \hbox{\rm in} \ \Omega, \\
\hspace{3.8cm} \oint_{\sigma_n} {\bf v} \cdot d{\bf s} = 0  \ \hbox{\rm for each} \ n = 1,\dots,g\} \, .
\end{array}
\end{equation}
Note that this space can be written as ${\mathcal V}_0 = \grad [H^1(\Omega)]$: in fact, the inclusion $\grad [H^1(\Omega)] \subset {\mathcal V}_0$ is obvious, while the inclusion ${\mathcal V}_0 \subset \grad [H^1(\Omega)]$ is the classical result concerning scalar potentials (see, e.g., Cantarella, DeTurck and Gluck~\cite{CDG02}, Alonso Rodr{\'\i}guez and Valli~\cite{AV14}).
The vector field ${\bf V}$ is thus a solution to 
\begin{equation}\label{eq:varpb2}
\begin{array}{ll}
{\bf V} \in {\mathcal V}_0 \ : \ \int_\Omega\boldsymbol \mu {\bf V}  \cdot \grad {\varphi} = -\int_\Omega {g} \, {\varphi}- \int_\Omega\boldsymbol \mu {\bf u}^\ast \cdot \grad {\varphi}\\[1ex]
\hspace{6cm} + \int_{\partial \Omega} {b} \, {\varphi}  \ \ \ \ \forall \ \varphi \in H^1(\Omega) \, .
\end{array}
\end{equation}
It is easy to see that ${\bf V}$ is indeed the {\sl unique} solution of that problem: in fact, assuming $g= b=0$, ${\bf u}^\ast = {\bf 0}$, and taking $\varphi$ such that $\grad \varphi = {\bf V}$, it follows at once $\int_\Omega \boldsymbol \mu {\bf V} \cdot {\bf V} = 0$, hence ${\bf V}= {\bf 0}$.

\subsection{Finite element approximation}\label{FEMapprox2}

The finite element approximation of problem (\ref{eq:cd_nor}) has two steps. The first one is finding a finite element potential ${\bf u}^\ast_h \in N_h$ such that
\begin{equation}\label{eq:u*h2}
\left\{\begin{array}{ll}
\curl {\bf u}^\ast_h  = {\bf J}_h &\hbox{\rm in} \ \Omega\\
\oint_{\sigma_n} {\bf u}_h^\ast \cdot d{\bf s} =\beta_n &\hbox{\rm for each} \ n = 1,\dots,g \, ,
\end{array}
\right.
\end{equation}
where ${\bf J}_h \in RT_h$ is the Raviart--Thomas interpolant $I^{RT}_h {\bf J}$ of ${\bf J}$ (we therefore assume that ${\bf J}$ is so regular that its interpolant $I^{RT}_h {\bf J}$ is defined; for instance, as already recalled, it is enough to assume ${\bf J} \in (H^{\frac{1}{2} + \delta}(\Omega))^3$, $\delta >0$: see Monk~\cite[Lemma 5.15]{Monk}). In Alonso Rodr{\'\i}guez and Valli~\cite{AV14} an efficient algorithm for computing ${\bf u}^\ast_h$ is described.

The second step is related to the numerical approximation of problem (\ref{eq:varpb2}). It is quite easy to find a finite element subspace of  ${\mathcal V}_{0}$ and a suitable finite element basis. The natural choice is clearly
\begin{equation}\label{eq:space_h2}
\begin{array}{ll}
{\mathcal V}_{0,h} = \{{\bf v}_h \in N_h \, | \, \curl{\bf v}_h = {\bf 0} \ \hbox{\rm in} \ \Omega, \\
\hspace{3.5cm} \oint_{\sigma_n} {\bf v}_h \cdot d{\bf s} = 0  \ \hbox{\rm for each} \ n = 1,\dots,g\} \, ,
\end{array}
\end{equation}
which can be rewritten as ${\mathcal V}_{0,h} = \grad L_h$. Since the dimension of this space is $n_v - 1$, a finite element basis is determined by taking $\grad \psi_{h,i}$, $i=1,\ldots,n_v-1$, $\psi_{h,i}$ being the basis functions of the finite element space $L_h$.

The finite element approximation of (\ref{eq:varpb2}) is easily obtained:
\begin{equation}\label{eq:varpb_fin2}
\begin{array}{ll}
{\bf V}_h \in {\mathcal V}_{0,h} \ : \ \ \int_\Omega\boldsymbol \mu {\bf V}_h  \cdot \grad {\varphi}_h = -\int_\Omega {g} \, {\varphi_h}- \int_\Omega\boldsymbol \mu {\bf u}_h^\ast \cdot \grad {\varphi}\\[1ex]
\hspace{6cm} + \int_{\partial \Omega} {b} \, {\varphi}_h  \ \ \ \ \forall \ \varphi_h \in L_h^\ast \, \, ,
\end{array}
\end{equation}
where 
\begin{equation}\label{FEspace_red2}
L^\ast_h = \hbox{\rm span} \{\psi_{h,i}\}_{i=1}^{n_v - 1} = \{\varphi_h \in L_h \, | \, \varphi_h(v_{n_v})=0\} \, .
\end{equation}
The corresponding algebraic problem is a square linear system of dimension $n_v - 1$, and it is uniquely solvable. In fact,  since ${\mathcal V}_{0,h} = \grad  L_h^\ast$, we can choose ${\varphi}_h^\ast \in L_h^\ast$ such that $\grad {\varphi}_h^\ast = {\bf V}_h$; from $(\ref{eq:varpb_fin2})$ we find at once ${\bf V}_h = {\bf 0}$, provided that $g=b = 0$, ${\bf u}^\ast_h={\bf 0}$.

The convergence of this finite element scheme is easily proved by following the arguments previously presented.

\begin{theo} Let ${\bf V} \in {\mathcal V}_0$ and ${\bf V}_h\in {\mathcal V}_{0,h}$ be the solutions of problem (\ref{eq:varpb2}) and (\ref{eq:varpb_fin2}), respectively. Set ${\bf u} = {\bf V} + {\bf u}^\ast$ and ${\bf u}_h= {\bf V}_h + {\bf u}_h^\ast$, where ${\bf u}^\ast \in H(\curl;\Omega)$ and ${\bf u}_h^\ast \in N_h$ are solutions to problem (\ref{eq:u*2}) and (\ref{eq:u*h2}), respectively. Assume that ${\bf u}$ and ${\bf J}$ are regular enough, so that the interpolants $I^{N}_h {\bf u}$ and $I^{RT}_h {\bf J}$ are defined. Then the following error estimate holds
\begin{equation}\label{eq:error2}
\|{\bf u}-{\bf u}_h\|_{H(\scriptsize \curl;\Omega)} \le c_0 (\|{\bf u} - I^{N}_h {\bf u}\|_{L^2(\Omega)} + \|{\bf J} - I^{RT}_h {\bf J}\|_{L^2(\Omega)}) \, .
\end{equation}
\end{theo}
{\bf Proof.} Since $L_h^\ast \subset H^1(\Omega)$, we can choose $\varphi= \varphi_h \in L_h^\ast$ in (\ref{eq:varpb2}). By subtracting (\ref{eq:varpb_fin2}) from (\ref{eq:varpb2}) we end up with
$$
\int_\Omega \boldsymbol\mu [({\bf V} + {\bf u}^\ast) - ({\bf V}_h + {\bf u}_h^\ast)] \cdot \grad \varphi_h = 0 \qquad \forall \ \varphi_h \in L_h^\ast \, ,
$$
namely, 
\begin{equation}\label{eq:consist2}
\int_\Omega \boldsymbol\mu ({\bf u} -  {\bf u}_h)\cdot \grad \varphi_h = 0 \qquad \forall \ \varphi_h \in L_h^\ast\, .
\end{equation}
Then, since ${\mathcal V}_{0,h} = \grad L_h^\ast$ and thus ${\bf V}_h = \grad \varphi_h^\ast$ for a suitable $\varphi^\ast_h \in L^\ast_h$, from (\ref{eq:consist2}) we find
$$
\begin{array}{ll}
c_1 \|{\bf u}-{\bf u}_h\|_{L^2(\Omega)}^2 \!\!\!&\le \int_\Omega \boldsymbol \mu ({\bf u}-{\bf u}_h) \cdot ({\bf u}-{\bf u}_h)\\ [1ex]
&=  \int_\Omega \boldsymbol \mu ({\bf u}-{\bf u}_h) \cdot ({\bf u} - {\bf V}_h - {\bf u}_h^\ast) \\[1ex]
&=  \int_\Omega \boldsymbol \mu ({\bf u}-{\bf u}_h) \cdot ({\bf u} - \grad\varphi_h^\ast - {\bf u}_h^\ast) \\ [1ex]
&=  \int_\Omega \boldsymbol \mu ({\bf u}-{\bf u}_h) \cdot ({\bf u} - \grad\varphi_h - {\bf u}_h^\ast) \\[1ex]
&\le c_2 \|{\bf u}-{\bf u}_h\|_{L^2(\Omega)} \|{\bf u} - {\boldsymbol\Psi}_h  - {\bf u}_h^\ast\|_{L^2(\Omega)} \quad \forall \ {\boldsymbol\Psi}_h \in {\mathcal V}_{0,h} \, .
\end{array}
$$
We can choose ${\boldsymbol\Psi}_h = (I^{N}_h {\bf u}- {\bf u}_h^\ast) \in {\mathcal V}_{0,h}$; in fact, $\curl (I^{N}_h {\bf u}) = I^{RT}_h(\curl {\bf u}) = I^{RT}_h {\bf J}={\bf J}_h$ and $\oint_{\sigma_n} I^{N}_h {\bf u}  \cdot d{\bf s} =\oint_{\sigma_n} {\bf u}  \cdot d{\bf s} = \beta_n$ for each $n=1,\ldots,g$.
Then we find at once $\|{\bf u}-{\bf u}_h\|_{L^2(\Omega)} \le \frac{c_2}{c_1} \|{\bf u} - I^{N}_h {\bf u}\|_{L^2(\Omega)}$.
Moreover, $\curl({\bf u}-{\bf u}_h) = {\bf J} - {\bf J}_h = {\bf J}-I^{RT}_h {\bf J}$, and the thesis follows. \qed

\vspace{.12cm}

Note that sufficient conditions for defining the interpolants of ${\bf u}$ and ${\bf J}=\curl {\bf u}$ is that they both belong to $(H^{\frac{1}{2} + \delta}(\Omega))^3$, $\delta >0$ (see Monk~\cite[Lemma 5.15, Theor.\ 5.41]{Monk}). This is for instance satisfied if $\boldsymbol \mu$ is a scalar Lipschitz function in $\overline \Omega$ and $b \in H^\gamma(\Omega)$, $\gamma > 0$ (see Alonso and Valli~\cite{AV99}).

\subsection{The algebraic problem}

The solution ${\bf V}_h \in {\mathcal V}_{0,h}$ is given by ${\bf V}_h = \sum_{i=1}^{n_v -1} V_i \grad\psi_{h,i}$. Hence the finite dimensional problem (\ref{eq:varpb_fin2}) can be rewritten as
\begin{equation}\label{eq:matrix2}
\begin{array}{ll}
\displaystyle\sum_{i=1}^{n_v-1} V_i \int_\Omega \boldsymbol{\mu} \grad \psi_{h,i} \cdot \grad \psi_{h,j} = - \int_\Omega g \, \psi_{h,j} - \int_\Omega\boldsymbol \mu {\bf u}^\ast_h \cdot \grad \psi_{h,j}\\
\hspace{6cm} + \displaystyle\int_{\partial \Omega} b \,  \psi_{h,j}  \, ,
\end{array}
\end{equation}
for each  $j=1,\ldots,n_v-1$. 

We have at once
\begin{theo}
The matrix ${\bf K}^\ast$ with entries 
$$
K^\ast_{ji} = \int_\Omega \boldsymbol{\eta} \grad \varphi_{h,i} \cdot \grad \varphi_{h,j}
$$
is symmetric and positive definite.
\end{theo}


\bibliographystyle{elsarticle-num}
\bibliography{curl-div}

\end{document}